\documentclass{article}
\usepackage{yfonts,calligra}
\usepackage{irule,amsmath,amssymb,amsthm}
\usepackage{graphicx}
\usepackage[all]{xy}
 
\newtheorem{theorem}{Theorem}
\newtheorem{lemma}{Lemma}
\newtheorem{corollary}{Corollary}
\newtheorem{definition}{Definition}
\def\<{\left <}
\def\>{\right >}
\def\({\left (}
\def\){\right )}
\def\bra {\left [}
\def\ket{\right ]}
\def\cbra{\left \{}
\def\cket{\right \}}

\def\skep{\hskip2em}

\DeclareSymbolFont{AMSb}{U}{msb}{m}{n}
\DeclareMathSymbol{\N}{\mathbin}{AMSb}{"4E}
\DeclareMathSymbol{\Z}{\mathbin}{AMSb}{"5A}
\DeclareMathSymbol{\R}{\mathbin}{AMSb}{"52}
\DeclareMathSymbol{\Q}{\mathbin}{AMSb}{"51}
\DeclareMathSymbol{\I}{\mathbin}{AMSb}{"49}
\DeclareMathSymbol{\C}{\mathbin}{AMSb}{"43}
\begin{document}

\title{Dynamic Topological Completeness for $\R^2$}
\author{David Fern\'andez-Duque\\
University of Barcelona\\
{\tt fernandez-duque@ub.edu.}
   }
   \date{}
   
   \maketitle

\begin{abstract}

Dynamic topological logic ({\bf DTL}) combines topological and temporal modalities to express asymptotic properties of dynamic systems on topological spaces. A dynamic topological model is a triple $\left <X,f,V\right >$, where $X$ is a topological space, $f:X\to X$ a continuous function and $V$ a truth valuation assigning subsets of $X$ to propositional variables. Valid formulas are those that are true in every model, independently of $X$ or $f$. A natural problem that arises is to identify the logics obtained on familiar spaces, such as $\R^n$. It \cite{slav} it was shown that any satisfiable formula could be satisfied in some $\R^n$ for $n$ large enough, but the question of how the logic varies with $n$ remained open.

In this paper we prove that any fragment of ${\bf DTL}$ that is complete for locally finite Kripke frames is complete for $\R^2$. This includes ${\bf DTL}^\bigcirc$; it also includes some larger fragments, such as ${\bf DTL_1}$, where ``henceforth" may not appear in the scope of a topological operator. We show that satisfiability of any formula of our language in a locally finite Kripke frame implies satisfiability in $\R^2$ by constructing continuous, open maps from the plane into arbitrary locally finite Kripke frames, which give us a type of bisimulation. We also show that the results cannot be extended to arbitrary formulas of ${\bf DTL}$ by exhibiting a formula which is valid in $\R^2$ but not in arbitrary topological spaces.

\end{abstract} 

\section{Introduction}

The system ${\bf S4}$ is a modal logic that has a natural topological interpretation. It uses the propositional language, adding the single modality $\Box$. For axioms, we have all tautologies, along with
$$\Box\alpha\to\alpha,$$
$$\Box\alpha\to\Box\Box\alpha$$
and
$$\Box(\alpha\to\beta)\to(\Box\alpha\to\Box\beta).$$
The inference rules are
$$\infer{\Box\alpha}{\alpha}$$
and modus ponens. For a more extensive treatment of ${\bf S4}$ (and other modal logics), see, for example, \cite{black}.

We then interpret $\Box$ as ``the topological interior of"; if $\left <X,{\bf T}\right>$ is a topological space and $V$ a valuation assigning a subset of $X$ to each propositional variable, $x\in X$ satisfies a formula $\Box\varphi$ if and only if a neighborhood of $x$ satisfies $\varphi$. This interpretation can be traced back to Tarski and others in the 1930s. In fact, Tarski proved in \cite{tarski} that ${\bf S4}$ is complete for $\R$ under these semantics; a streamlined proof can be found in \cite{mints}.\\

Another class of topological spaces for which ${\bf S4}$ is complete is that of finite Kripke frames (\cite{black}; throughout this paper, all Kripke frames will be transitive and reflexive). A {\em Kripke frame} is a pair
$${\bf K}=\<W,R\>,$$
where $W$ is a set and $R$ a preorder on $W$ (that is, a transitive, reflexive relation). If $Rwv$, we will often say that $v$ is a {\em successor} of $w$. Such a relation naturally induces a topology, by letting $U\subset W$ be open if and only if, whenever $w\in U$, all successors of $w$ lie in $U$. These semantics were introduced by Saul Kripke in 1959. Historically, topological semantics precede Kripke semantics, but the latter are now very standard in modal logic.

Closely related to modal logics are temporal logics, introduced in the 1960s by Arthur Prior and other authors. These systems model the dependence of truth on time. In the context of this paper `time' is a discrete secuence of states, each one indexed by a natural number, and modalities refer to future events; $\bigcirc$ stands for `next', so that $\bigcirc p$ is true on state $n$ if and only if it $p$ is true in state $n+1$, while $\ast$ stands for `henceforth'; $\ast p$ is true in state $n$ if and only if $p$ holds in all states $n+k$ for $k\in\N$.

We are interested in extensions of {\bf S4} which also use temporal operators. Such combined systems were introduced in the late 1990s (\cite{arte,kmints}) with the purpose of reasoning about dynamic systems on topological spaces. A {\em dynamic topological model} is a pair $\left<\left<X,{\bf T}\right>,f\right>$, where $\left<X,{\bf T}\right>$ is a topological space and $f:X\to X$ a function. Specifically, we are interested in the case where $f$ is continuous, and by {\em Dynamic Topological Logic} (${\bf DTL}$) we will refer exclusively to logics interpreted on such systems.

The language of {\bf DTL} uses the language of propositional logic along with the three modal operators discussed above; $\Box,\bigcirc$ and $\ast$.

To each propositional variable $p$ we assign an arbitrary set $V(p)\subset X$. Boolean connectives are interpreted as usual, while $V(\Box\alpha)=V(\alpha)^\circ$ (the topological interior), $V(\bigcirc\alpha)=f^{-1}(V(\alpha))$ and
$$V(\ast\alpha)=\bigcap_{n=0}^\infty f^{-n}(V(\alpha)).$$
We say that
$$x\models \varphi\Leftrightarrow x\in V(\varphi).$$

A formula $\varphi$ is valid if it is valid in every dynamic topological system (that is, $V(\varphi)=X$). The set of valid formulas of ${\bf DTL}$ is undecidable (\cite{konev}). However, some of its fragments turn out to be rather manageable. We are especially interested in fragments which are complete for Kripke frames.

One such fragment is ${\bf DTL^{\bigcirc}}$, the $\ast$-free fragment of ${\bf DTL}$. A sound and complete axiomatization for it was introduced in \cite{kremer} (there called {\bf S4C}). It includes all axioms and rules of ${\bf S4}$, along with
$$\bigcirc(\alpha\odot\beta)\leftrightarrow\bigcirc\alpha\odot\bigcirc\beta$$
for any Boolean connective $\odot$,
$$\bigcirc\Box\alpha\to\Box\bigcirc\alpha,$$
which expresses continuity, and the inference rule
$$\infer{\bigcirc\alpha}{\alpha}.$$

Another fragment of interest is ${\bf DTL_1}$, which properly includes ${\bf DTL}^\bigcirc$. ${\bf DTL_1}$ uses all three operators but $\ast$ may never appear in the scope of $\Box$. This fragment is undecidable but recursively enumerable ($\cite{konev}$).

It does have the following property, which we will rely heavily on:

\begin{theorem}\label{aleks}
A formula of ${\bf DTL_1}$ is satisfiable if and only if it can be satisfied in a locally finite Kripke frame.
\end{theorem}

``Locally finite" here means that each world has finitely many successors.

\proof
A proof can be found in \cite{konev}.
\endproof

There is no generalization of Tarski's theorem for ${\bf DTL}^{\bigcirc}$ on the real line, as was observed by P. Kremer in \cite{kremer} and S. Slavnov in \cite{slav2}. This raises the question of completeness of the system on $\R^n$ for $n>1$. In \cite{slav}, S. Slavnov shows that it is complete for $\amalg _{n=1}^\infty \R^n$, but the question of completeness for any fixed $n$ remained open. This paper gives a positive answer; $n=2$ is, in fact, sufficient.

Our strategy is based on Slavnov's. Essentially, one exploits the completeness of ${\bf DTL}^\bigcirc$ for finite Kripke frames $K=\left<W,R,g\right >$ and generates continuous, open mappings $$\omega:\R^2\to W.$$ Such mappings preserve the truth of formulas of ${\bf DTL}$; see \cite{kremer} or \cite{slav}.

The main difference with \cite{slav} is in the function
$$f:\R^2\to\R^2$$
that we use. The reason ${\bf DTL}^{\bigcirc}$ is not complete for $\R$ is that there is too little room to
navigate between a constant map and an open map (see the counterexample given in \cite{kremer}); thus our construction
must be such that $f$ is nowhere open, and yet never collapses an open ball
to a one-dimensional space, even after many iterations.
It is not at first evident how one is to define such a function,
and most of this paper is dedicated to setting up the necessary tools.

Our goal is to prove that any formula of ${\bf DTL}$ that can be satisfied in a locally finite Kripke frame can be satisfied in $\R^2$. This gives us as immediate corollaries that satisfiability in ${\bf DTL}^\bigcirc$ and ${\bf DTL_1}$ is equivalent to satisfiability in $\R^2$. In Section 2, we define the notion of a topological bisimulation, which is essentially a continuous, open map which preserves valuations of {\bf DTL} formulas. In Section 3 we define open ball trees and give a proof of Tarski's theorem for the plane. The proof we give is very similar to previously existing versions but is adjusted so it can be extended to work with temporal operators. Section 4 deals with dynamic Kripke frames; eventually we wish to construct a topological bisimulation from $\R^2$ to a dynamic Kripke frame. Section 5 defines segment trees, which like open ball trees are auxiliary tools meant to encode information from a Kripke frame geometrically on the plane. Section 6 defines dynamic sequences of open ball trees and proves our main theorem. Finally, in Section 7 we define a formula of {\bf DTL} which is valid in all locally connected topological spaces but not in arbitrary spaces; this shows that our results are sharp in the sense that they cannot be extended to arbitrary formulas of {\bf DTL}. Some tools from analysis are needed, but all are fairly standard; Appendix A gives an overview of the results and conventions used throughout the paper.

\section {Topological bisimulation}

The following theorem can be found in $\cite{kremer}$ or $\cite{slav}$ and will be essential for obtaining completeness results in $\R^2$:
\begin{theorem}\label{prel}
Let $$\left< X,f_X,V_X \right >$$
and $$\left< Y,f_Y,V_Y \right >$$
be dynamic topological models and
$$\gamma:Z\subset X\to Y$$
a continuous and open map such that
$$f_Y\gamma=\gamma f_X$$
and
$$V_X=\gamma^{-1}V_Y.$$
Then, for any ${\bf x}\in Z^\circ$ and an arbitrary ${\bf DTL}$ formula $\varphi$,
$${\bf x}\models_X \varphi\Leftrightarrow \gamma({\bf x})\models_Y\varphi.$$
\end{theorem}
\proof
The proof proceeds by induction on the depth of $\varphi$. The induction step for Boolean connectives works trivially. For $\varphi=\Box\psi$, note that $x\models_X \Box\psi$ if and only if there is a neighborhood $U_x$ of $x$ such that $x'\models_X \psi$ for all $x'\in U_x$. But using our induction hypothesis and the fact that $\gamma$ is continuous and open, this is equivalent to there being a neighborhood $U_{\gamma(x)}$ of $\gamma(x)$ such that all points of $U_{\gamma(x)}$ satisfy $\psi$.

If $\varphi=\bigcirc\psi$, we have that $x\models_X \varphi$ if and only if $f_X(x)\models_X\psi$. This is equivalent by induction hypothesis to $\gamma(f_X(x))\models_Y\psi$, and this in turn to $f_Y(\gamma(x))\models_Y\psi,$ or $\gamma(x)\models_Y\bigcirc\psi.$

Finally, for $\varphi=\ast\psi$,
$$
\begin{array}{lcl}
\gamma^{-1}(V_Y(\ast\psi))&=&\gamma^{-1}\left (\bigcap_{n=0}^\infty f_Y^{-n}\left (V_Y(\psi)\right )\right )\\\\
&=&\bigcap_{n=0}^\infty \gamma^{-1}\left (f_Y^{-n}\left (V_Y(\psi)\right )\right )\\\\
&=&\bigcap_{n=0}^\infty f^{-n}\left (\gamma^{-1}\left (V_Y(\psi)\right )\right )\\\\
\text{\skep\small\bf IH}&=&\bigcap_{n=0}^\infty f^{-n}\left (V_X(\psi)\right )\\\\
&=&V_X(\ast\psi).
\end{array}
$$
\endproof
\section {Open ball trees in $\R^2$ and completeness of ${\bf S4}$}

Throughout this section, $K=\left<W,R\right>$ is a finite, transitive, reflexive Kripke frame (that is, $W$ is a finite set and $R$ a transitive, reflexive relation on $W$). It is a well-known result that ${\bf S4}$ is complete for the class of such frames, and we can exploit this to prove completeness for $\R^2$ (in fact the same proof works for $\R$).
\begin{figure}[htp]
\begin{center}
\scalebox
{0.9}
{
\includegraphics[viewport=195 285 600 330]{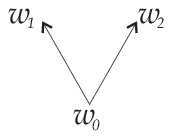}
}
\end{center}
\caption
{
A simple Kripke frame.
}
\label{kripke}
\end{figure}\\

Note that $W$ naturally acquires a topology; open sets are sets that are upward closed under $R$. Unfortunately such a space is not Hausdorff, and as such a converging sequence does not necessarily have a unique limit. We remedy this, following Mints, by assigning to each $R$-monotone sequence $\{w_n\}_{n=0}^\infty$ a single limit $w\in\{w_n\}_{n=0}^\infty$ in such a way that $Rw_nw$ and, for $n$ large enough, $Rww_n$. Because $W$ is locally finite it is clear that such a $w$ can always be selected; denote it by $\lim_{n\to\infty}w_n$.

Our strategy for proving completeness of ${\bf S4}$ for $\R^2$
will be to define an open and continuous map from $\R^2$ to our Kripke frame $K$ and apply Theorem \ref{prel}. The following construction will be useful in defining such a map:\\\\
\begin{definition}[Open ball tree]
An {\em open ball tree} is a set $\Theta$ whose elements are open balls in $\R^2$ satisfying
the following three properties:\\
\\ {\bf (obt1)} $\Theta$ has a maximum element under inclusion, which we will denote $root(\Theta)$;\\
\\ {\bf (obt2)} if $B_1,B_2\in \Theta$ and $\partial B_1\cap \partial B_2\not =\emptyset $,
then $B_1=B_2$;\\
\\ {\bf (obt3)} given $\varepsilon>0$, there are only finitely many elements of $\Theta$ with diameter $\geq\varepsilon$.
\end{definition}

It is not hard to see that the elements of an open ball tree, ordered by inclusion, form a tree, with the largest element as
root.
\begin{figure}[htp]
\begin{center}
\scalebox
{0.9}
{
\includegraphics[viewport=195 280 600 350]{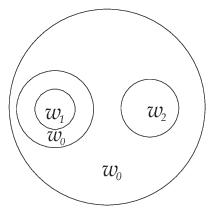}
}
\end{center}
\caption{
A typical open ball tree. We have labeled each ball by a world of the Kripke frame from Figure
\ref{kripke}; this labeling represents a continuous map $\omega$.
}
\label{typical}
\end{figure}

Because of condition {\bf (obt2)}, if $B$ and $B'$ contain a point ${\bf x},$ then either $B\subset B'$
or $B'\subset B$. Thus the open balls containing ${\bf x}$ are totally ordered. Also note that, because of {\bf (obt3)}, given any open ball $B$ containing ${\bf x}$, there are only finitely many balls $B'\in \Theta$ such that $B\subset B'$ and ${\bf x}\in B'$. Hence all open balls of ${\Theta}$ containing ${\bf x}$ form a sequence, 
$$root(\Theta)=B_0\supset B_1\supset B_2\supset...$$

We will call the elements of this sequence $\theta_n({\bf x})$. Note that this sequence may be finite or infinite; if it is finite, we will denote the minimum element by $\theta({\bf x})$.

Otherwise, using {\bf (obt3)} we know that
$$\lim_{n\to\infty}diam\(\theta_n({\bf x})\)=0,$$
so $\{{\bf x}\}=\bigcap_{n=0}^\infty \theta_n({\bf x})$.
In this case we will say ${\bf x}$ is a {\em limit point} of $\Theta$, and denote the set of all such points by
$lim(\Theta)$.

We will also use the notation
$$dom(\Theta)=root(\Theta)\setminus lim(\Theta).$$

A mapping $\omega : \Theta\to W$ is continuous if
$B_1\supset B_2\Rightarrow R\omega (B_1)\omega(B_2)$ (this is simply
continuity in the topology of preorders).
A continuous map $\omega$ induces a mapping $\dot \omega :root(\Theta)\to W$
defined as follows:\\\\
$$
\dot \omega({\bf x})=
\begin{cases}
\omega(\theta({\bf x}))&\text{if ${\bf x}\in dom(\Theta)$;}\\\\
\lim_{n\to\infty}\omega(\theta_n({\bf x}))&\text{if ${\bf x}\in lim(\Theta)$.}
\end{cases}
$$
\\\\
We then have the following lemma:
\begin{lemma}
If $\omega$ is continuous, then $\dot\omega$ is continuous as well.
\end{lemma}
\proof
Let ${\bf x}\in root(\Theta)$. First assume ${\bf x}\in dom(\Theta)$, and let
${\bf y}\in \theta({\bf x}).$ If ${\bf y}\in dom(\Theta)$, then
evidently $\theta({\bf x})\supset\theta({\bf y})$ so $R\omega({\bf x})\omega({\bf y})$.
Otherwise, since $\theta({\bf x})=\theta_n({\bf y})$ for some $n$, we have
$R\omega(\theta({\bf x}))\lim_{n\to\infty}\omega(\theta_n({\bf y}))$ (by definition of limits on Kripke frames).
In either case we conclude $R\dot\omega({\bf x})\dot\omega(\bf y)$, which means ${\bf y}$ is contained in the preimage
of any open set containing $\dot\omega({\bf x})$. Hence $\dot\omega$ is continuous at ${\bf x}$.

If ${\bf x}$ is a limit point, then pick $n$ large enough so that $R\dot\omega({\bf x})\omega(\theta_n({\bf x})).$
The same argument as before goes through with $\theta_n({\bf x})$ in place of $\theta({\bf x})$.
\endproof

Continuity of $\omega$ does not guarantee that it is open. For this we need a bit more work.
\begin{definition}[$\varepsilon$-grid]
Let $\varepsilon>0$, and $U\subset \R^2$ be an open set. The $\varepsilon$-grid on $U$
is the collection of connected components of $$U\cap\bigcup_{n,m\in{\bf Z}}\left \{\left(
\begin{array}{c}
x\\
y\\
\end{array}
\right) \in\R^2 : n\varepsilon<x<(n+1)\varepsilon \text{ and }m\varepsilon<y<(m+1)\varepsilon\right \}.$$
We will refer to it as $\Gamma_\varepsilon (U)$.
\end{definition}

Thus we chop the set $U$ into many small open pieces, throwing away their boundaries. These grids will be
essential to construct open mappings into our Kripke frame.

\begin{definition}[$\varepsilon$-dense]
Let $\Theta^0\subset\Theta^1$, $\omega$ a continuous map defined on $\Theta^1$ and $\varepsilon>0$.

We say $\Theta^1$ is $\varepsilon$-dense
in $\Theta^0$ if the following conditions hold:\\\\
{\bf (dns1)} $\Theta^1$ is finite;\\\\
{\bf (dns2)} if ${\bf x}\in root(\Theta^0)$ and $R\omega(\theta^1({\bf x}))w$ there exists $B\in \Theta^{1}$ such that
$\omega(B)=w$ and $d({\bf x},B)<\varepsilon$;\\\\
{\bf (dns3)} if $B\in \Theta^{1}\setminus\Theta^0$ then $diam (B)< \varepsilon$.
\end{definition}
\begin{definition}[dense sequence]
A pair $(\left \{\Theta^m\}_{m=0}^\infty, \omega\right )$ where $\{\Theta^m\}_{m=0}^\infty$ is a
sequence of open ball trees such that $\Theta^m\subset\Theta^{m+1}$
and
$$\omega:\Theta=\bigcup_{m=0}^\infty \Theta^m\to W$$
is a continuous mapping is {\em dense} if
$\Theta^{m+1}$ is $\left (\frac 12\right )^m$-dense in $\Theta^m$.
\end{definition}

For an open ball tree
$\Theta$, we will use the notation $$\Theta^\circ=root(\Theta)\setminus\bigcup_{B\in\Theta}\partial B$$
and refer to $\Theta^\circ$ as the {\em interior} of $\Theta$. As long as $\Theta$ is finite, $\Theta^\circ$ will be an open set.

\begin{figure}[htp]
\begin{center}
\scalebox
{0.9}
{
\includegraphics[viewport=195 280 600 350]{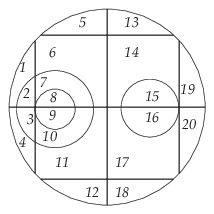}
}
\end{center}
\caption{
A grid on the interior of the open ball tree from Figure \ref{typical}. We have numbered each element of the grid.
}
\label{grid}
\end{figure}

\begin{definition}[$\varepsilon$-refinement]
Let $\varepsilon>0$, $\Theta$ a finite open ball tree and $\omega:\Theta\to W$ a continuous map.

Let $C_1,...,C_I$ be the open, connected sets that make up the $\varepsilon$-grid $\Gamma_\varepsilon\(\Theta^\circ\)$.
Note that $\dot\omega$ is a constant on each component of $\Theta^\circ$ and hence
constant on each $C_i$. We will denote this constant by $\omega(C_i)$.

Then, $(\check \Theta,\check\omega)$ is an {\em $\varepsilon$-refinement} of $(\Theta,\omega)$ if\\\\
{\bf (ref1)} $\check\Theta$ is finite;\\\\
{\bf (ref2)} $\Theta\subset\check\Theta;$\\\\
{\bf (ref3)} if $B\in \check\Theta\setminus\Theta$ then $\overline B\subset C_i$ for some $i$;\\\\
{\bf (ref4)} if $B_0,B_1\in\check\Theta\setminus\Theta$ and $B_0\not=B_1$, $\overline{B_0}\cap\overline{B_1}=\emptyset;$\\\\
{\bf (ref5)} $\check\omega$ is continuous and $\check\omega|_{\Theta}=\omega$;\\\\
{\bf (ref6)} Whenever
$R\omega(C_i)w$, there exists $B\in \check\Theta$ contained in $C_i$ such that $\check\omega(B)=w.$
\end{definition}

An $\varepsilon$-refinement always exists and is not hard to construct.
Let $w_{i,0},...,w_{i,N_i}$ be all the successors
of $\omega(C_i)$. Since $C_i$ is open, we can pick out $N_i$ open balls $B_{i,n}$ whose closures
are contained in $C_i$ and pairwise disjoint. Then set $\check\omega(B_{i,n})=w_{i,n}$ and
$\check\Theta=\Theta\cup\bigcup_{i,n}\{B_{i,n}\}$.\\
\begin{figure}[htp]
\begin{center}
\scalebox
{0.9}
{
\includegraphics[viewport=200 260 600 360]{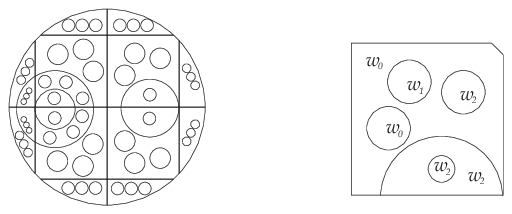}
}
\end{center}
\caption{
A refinement using the grid from Figure \ref{grid}. On the right hand side we zoom in
on to a portion of the grid to show the mapping $\check\omega$.
}
\end{figure}
\begin{lemma}
Let $\{\Theta^m\}$ be a sequence of open ball trees and $\omega:\bigcup_{m=0}^\infty\Theta^m\to W$ a continuous
map such that $(\Theta^{m+1},\omega|_{\Theta^{m+1}})$ is a $\left (\frac{1}{2}\right )^{m+1}$-refinement
of $(\Theta^{m+1},\omega|_{\Theta^{m+1}})$. Then, $(\{\Theta^m\},\omega)$ is dense.
\end{lemma}
\proof
We will prove that $\Theta^{m+1}$ is $\left (\frac 12\right )^m$-dense in $\Theta^m$.\\\\
{\bf (dns1)},{\bf (dns2)}: Obvious by definition of refinements.\\\\
{\bf (dns2)}: Let $C_1,C_2,...,C_I$ be the elements of the $\(\frac 12 \) ^{m+1}$-grid of $(\Theta^m)^\circ$,
${\bf x}\in root(\Theta^0)$ and $R\omega(\Theta^m({\bf x}))w$. Pick out an
index $i$ such that ${\bf x}\in \overline C_i$.
Then $C_i\subset B_{\left (\frac 12 \right )^m}({\bf x}),$
and by {\bf (ref6)} there is some $B\in \Theta^{m+1}$
contained in $C_i$ and such that $\omega(B)=w$.\\\\
{\bf (dns3)}: Evident once we consider that all the balls added at stage $m+1$ are
contained in one component of the
$\left (\frac 12 \right )^{m+1}$-grid.
\endproof

It is easy to see that if $(\{\Theta^m\}_{m=0}^\infty,\omega)$ is dense,
$\Theta=\bigcup_{m=0}^\infty\Theta^m$ is
an open ball tree and $\omega$ a continuous map. One can also observe that if
$\omega(B)=w$, then $B$ contains points
such that $\dot\omega({\bf x})=w$ (for example, take a point that is on
the boundary of any ball contained in $B$). This
last consideration is enough to prove that $\dot\omega$ is open:
\begin{lemma}\label{old}
The map $\dot\omega$ thus defined is continuous and open.
\end{lemma}
\proof We have already proven that such maps are continuous. Now, to see that
it is open, take ${\bf x}\in B_1(0)$,
$\varepsilon>0$ and $w\in W$ such that $R\dot\omega({\bf x})w$.
Let $m$ be such that $\left (\frac 12 \right )^m<\varepsilon$.
By {\bf (dns2)} there is $B\in \Theta^m$ such that $B\subset B_{\left (\frac 12 \right )^m}({\bf x})$
and $\omega(B)=w$.

Since
$\dot\omega({\bf y})=w$ for some ${\bf y}\in B$, we have that
$$w\in \dot\omega \left (B_{\left (\frac 12 \right )^m}\left({\bf x}\right)\right ).$$ 
\endproof
\begin{theorem}
${\bf S4}$ is complete for $\R^2$.
\end{theorem}
\proof
Let $\varphi$ be a non-valid formula of ${\bf S4}$. Let $M=\left <W,R,V,w_0\right >$ be a finite countermodel
for $\varphi$ with limits.

Define $\Theta^0=\{B_1({\bf 0})\}$ (the unit ball centered at the origin) and $\omega(B_1({\bf 0}))=w_0$.

Construct a sequence $\{(\Theta^n,\omega|_{\Theta^n})\}_{n=0}^\infty$
such that $(\Theta^{n+1},\omega|_{\Theta^{n+1}})$ is a $\left (\frac{1}{2}\right )^{n+1}$
-refinement of $(\Theta^n,\omega|_{\Theta^n})$.

Then, $\omega$ defines a continuous and open mapping $\dot\omega:B_1({\bf 0})\to W$. Pick ${\bf x}\in B_1({\bf 0})$
such that $\dot\omega({\bf x})=w_0$.

By Theorem \ref{prel},
$$\<\R^2,\dot\omega^{-1}V,{\bf x}\>\not\models\varphi.$$
\endproof
\section {Dynamic Kripke frames}

${\bf DTL}^{\bigcirc}$ is complete for the class of finite dynamic Kripke models (\cite{arte}). These are tuples
$$\left< W,R,g,V \right >$$
where $W$ is a finite set, $R$ a transitive, reflexive relation on $W$ and $g$
an $R$-monotone function on $W$. $W$ naturally
acquires a topology; open sets are simply $R$-closed subsets of $W$, and it can be easily checked that monotonicity of $g$ is equivalent to continuity under this topology.

${\bf DTL_1}$ is not complete for the class of finite Kripke frames, but it is for the class of locally finite Kripke frames ($\cite{konev}$). These are Kripke frames where every world has finitely many successors. 
\begin{figure}[htp]
\begin{center}
\scalebox
{0.9}
{
\includegraphics[viewport=190 260 600 360]{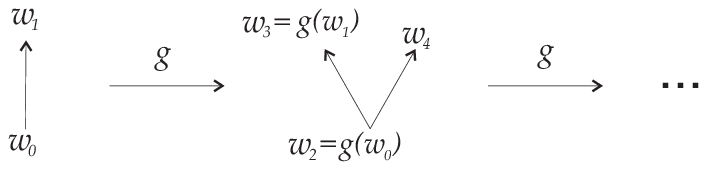}
}
\end{center}
\caption{
A portion of a dynamic Kripke frame.
}
\label{dynk}
\end{figure}

We would like the function $g$ to commute with limits in the usual way. The following theorem can essentially be found in \cite{arte} in a different presentation.
\begin{theorem}
Let $M=\left< W,R,g,V \right >$ be a locally finite Kripke model and $w_*\in W$ fixed. Then there exists a locally finite dynamic Kripke model
$$\tilde M=\left< \tilde W,\tilde R,\tilde g,\tilde V \right >$$
where limits are defined,
$$\tilde g\left (\lim_{n\to\infty}w_n\right )=\lim_{n\to\infty}\tilde g(w_n),$$
every world has finitely many successors and such that there exists $\tilde w_*\in\tilde W$ satisfying
$$w_*\models_M\varphi\Leftrightarrow \tilde w_*\models_{\tilde M}\varphi$$
for any ${\bf DTL}$ formula $\varphi$.
\end{theorem}
\proof
Start with a locally finite model
$$M=\left< W,R,g \right >$$
with arbitrary limits. It might not be possible to make limits commute with $g$ in the original model. Instead, following \cite{slav} we begin by ``stratifying" $W$; that is, we will construct a set
$$\tilde W=\bigcup_{n=0}^\infty W_n$$
such that the $W_n$ are disjoint and open.

Let
$$W_0=\left \{(w):Rw_*w\right \},$$
and define inductively
$$W_{n+1}=\left \{\left (w_0,...,w_{n},w_{n+1}\right ):\left (w_0,...,w_{n}\right )\in W_n\text{ and }Rg(w_n)w_{n+1}\right \}.$$
Define
$$\tilde R\left (w_0,...,w_{n}\right )\left (v_0,...,v_{m}\right )\Leftrightarrow n=m\text{ and }Rw_kv_k\text{ for all }k\leq n,$$
$$\tilde g\left (w_0,...,w_{n}\right )=\left (w_0,...,w_{n},g(w_n)\right )$$
and
$$\tilde V(p)=\left \{ \left (w_0,...,w_{n}\right ):w_n\in V(p)\right \}.$$

Finally, set
$$\lim_{k\to\infty}\left (w_{0,k},...,w_{n,k},g\left (w_{n,k}\right )\right )=\left (\lim_{k\to\infty}w_{0,k},...,\lim_{k\to\infty}w_{n,k},g\left (\lim_{k\to\infty}w_{n,k}\right )\right ).$$

It is not hard to see that
$$\tilde M=\left <\tilde W,\tilde R,\tilde g,\tilde V\right >$$
satisfies all the desired properties, with $\tilde w_*=(w_*).$
\endproof
\section{Segment trees in $\R^2$}

To prove completeness results for dynamic topological logics for $\R^2$ we will use, along with open ball trees,
a second auxiliary construction which we will call {\em segment trees}. These can be thought of
as a line segment in $\R^2$, out of which many smaller segments sprout. Out of each of the smaller
segments many more may sprout, so that they look like smaller copies of the whole tree, and this
can continue indefinitely.\\\\
\begin{definition}[segment tree]
A {\em segment tree} is a set $\Phi $ of vertical and
horizontal closed line segments in $\R^2$ (which we will refer to as {\em branches}) such that:\\\\
{\bf (seg1)} There is an element $root(\Phi)\in\Phi$ such that its endpoints do not intersect any other
element of $\Phi$;\\\\
{\bf (seg2)} If $b_0,b_1\in\Phi$ and $b_0\cap b_1\not=\emptyset$, then the intersection is a single point
which is an endpoint of $b_0$ or $b_1$, but not of both.\\\\
{\bf (seg3)} For any $b\not= root(\Phi)$, exactly one of its endpoints intersects another element $b'\in\Phi$.
We write $b\dashv b'$.\\\\
{\bf (seg4)} There is no infinite sequence such that
$$b_0\dashv b_1\dashv b_2\dashv...$$
\end{definition}
\begin{figure}[htp]
\begin{center}
\scalebox
{0.9}
{
\includegraphics[viewport=200 260 600 350]{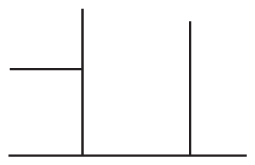}
}
\end{center}
\caption{
A segment tree. Compare this to the open ball tree of Figure \ref{typical}; both represent the same
tree.
}
\end{figure}

Segment trees are another way to represent trees geometrically in $\R^2$, where the partial
order is the transitive, reflexive closure of $\dashv$ (which we will denote $\dashv^*$).
In this way they are similar to open ball trees,
and when an open ball tree and a segment tree represent the same tree, one can find
a correspondence between open balls and branches which respects the partial order.

We will use the notation $im(\Phi)={\bigcup_{b\in \Phi}b}$.

We have a function $\phi:im(\Phi)\to\Phi$ assigning to each point the branch that it is on.
If ${\bf y}\in b_0\cap b_1$ and $b_1\dashv b_0$, set $\phi({\bf y})=b_0$.\\\\
\begin{definition}[tree maps]\label{treemaps}
Let $\Theta$ be an open ball tree and $\Phi$
a segment tree. Let $f:root(\Theta)\to\R^2$ be a continuous function.

We say that $f$ is a {\em tree map} 
between $\Theta$ and $\Phi$ if there is an isomorphism
$f_- :\Theta\to\Phi$ making the following diagram commute:\\
\begin{equation}
\xymatrix{
\Theta\ar@{-->}[r]^-{f_-}_-{\cong} & \Phi\\
dom(\Theta )\ar[u]^-{\theta}\ar[r]^-{f} & im(\Phi).\ar[u]^-{\phi}
}
\end{equation}

We will write $\xymatrix{ \Theta\ar[r]^-{f} & \Phi}.$
\end{definition}
 
Here {\em isomorphism} refers to an isomorphism as partially ordered sets.

Tree maps can be built whenever we have two
structures representing the same finite tree.
We do this by a procedure we will call {\em pinching}.\\\\
\begin{definition}[pinching]
Let $f:U\subset \R^2\to \R^2$ be a continuous function,
$B_0,B_1\subset U$ open balls such that $\overline B_1\subset B_0$ and $\overline B_0\subset U$;
let $b_0$, $b_1$ be branches such that $b_1\dashv b_0$. Assume $f(B_0)\subset b_0$. We say that a function
$\hat f$ is obtained from $f$ by {\em pinching} $B_1$  at $B_0$ onto $b_1$ if
\\\\
{\bf (pin1)} $\hat f$ is continuous on $U$;\\\\
{\bf (pin2)} $\hat f|_{U\setminus B_0 }=f|_{U\setminus B_0 };$\\\\
{\bf (pin3)} $\hat f(B_0\setminus B_1)\subset b_0;$\\\\
{\bf (pin4)} $\hat f(B_1)=b_1\setminus b_0$. \\\\
We will write $\hat f\succ _{B_0}^{B_1\to b_1}f$.
\end{definition}
\begin{lemma}\label{pinch}
With the notation of the previous definition assume $b_1$ has length less than $\frac{\varepsilon}{3}$ and $f(B_0)$ is contained in 
a ball of radius $\frac{\varepsilon}{3}$. Then there exists a function
$\hat f\succ _{B_0}^{B_1\to b_1}f$ such that $\parallel\hat f-f\parallel_\infty<\varepsilon$.
\end{lemma}
Recall that
$$\parallel f\parallel_\infty=\sup_{x\in dom(f)}\|f(x)\|.$$
\proof
Let ${\bf y}$ be the point of intersection of $b_0$ and $b_1$. First we find a function $\tilde f$ such
that $$\parallel \tilde f-f\parallel _\infty < \frac{2}{3}\varepsilon,$$
$\tilde f({\bf x})={\bf y}$ (where ${\bf x}$ is the center of $B_1$), and $\tilde f=f$ outside of $B_0$; such a function can be constructed by a small deformation of $f$, occurring within the ball of radius $\frac\varepsilon 3.$

We then build $\hat f$ by a similar deformation, as shown in Figure \ref{figp} below. Because $b_1$ has length less than $\frac\varepsilon 3$, we can ensure that
$$\parallel\hat f-f\parallel_\infty<\varepsilon,$$
by taking all deformations within the union of $b_1$ and the ball of radius $\frac\varepsilon 3$.
\begin{figure}[htp]
\begin{center}
\scalebox
{0.9}
{
\includegraphics[viewport=210 270 600 380]{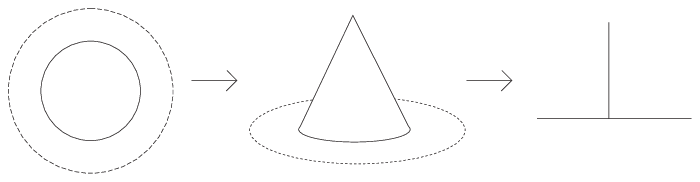}
}
\caption
{
Pinching can be thought of in two steps, where we first pull $B_1$ onto a cone and then collapse the respective
parts onto line segments.
}
\label{figp}
\end{center}
\end{figure}
\endproof

Iterating this process we can get tree maps from a finite
open ball tree to a finite segment tree. In fact we can iterate this infinitely many times
as long as the intermediate functions converge. Eventually, we will want to produce a tree map from an infinite open ball tree to an infinite segment tree.\\
\begin{figure}[htp]
\begin{center}
\scalebox
{0.9}
{
\includegraphics[viewport=200 200 600 400]{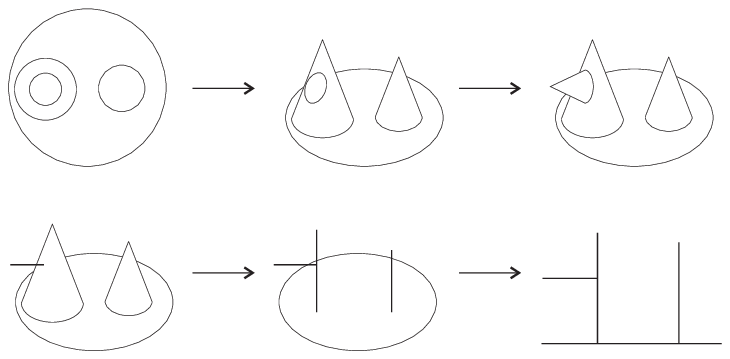}
}
\end{center}
\caption{
By iterating pinching, we can generate tree maps between more complex structures.
}
\end{figure}

There is a second way in which segment trees can be related to open ball trees.\\\\
\begin{definition}[embedded segment trees]\label{embedded}
{\em Let $\Phi$ be a segment tree and $\Theta$ be an open ball tree.
Suppose $im(\Phi)\subset dom(\Theta).$
We say $\Phi$ is {\em embedded} in $\Theta$ if there is a mapping $\iota : \Phi\to\Theta$
making the following diagram commute:
\begin{equation}
\xymatrix{
\Phi\ar@{-->}[r]^-{\iota} & \Theta\\
im(\Phi)\ar[u]^{\phi}\ar@{^{(}->}[r] & dom(\Theta)\ar[u]^-{\theta}
,}
\end{equation}
and such that whenever $\iota(b)\subset B\in\Theta$, there exists $b^*$ such that $\iota(b^*)=B$ and $b\dashv^*b^*$.

We will write $\Phi\hookrightarrow \Theta$.}
\end{definition}

\begin{figure}[htp]
\begin{center}
\scalebox
{0.9}
{
\includegraphics[viewport=200 260 600 370]{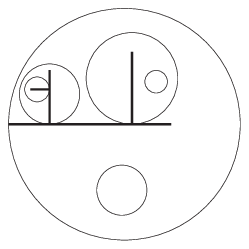}
}
\end{center}
\caption{
A segment tree embedded in an open ball tree; $\iota(b)$ is the smallest ball
containing $b$. Note that many of the balls are tangent to a segment;
this is necessary to construct an embedded segment tree. Note also that not all balls have a segment
corresponding to them.
}
\label{tangent}
\end{figure}

One can construct a segment tree embedded in an open ball tree as follows:\\\\
Start with an open ball $B_0$ and a line segment $b_0$ which is one of the radii of $B_0$.
To add a branch and an open ball, pick $B$ tangent to $b_0$ and small enough so that $\overline B\subset B_0.$
Then add the radius of $B$ that touches $b_0$ as a new branch $b$.

In further iterations, we can add other branches perpendicular to $b_0$, or we can also add branches onto the new
branches, contained in their respective open balls.

In fact, we have the following:
\begin{lemma}\label{addbranch}
With notation as above, let $\Phi\hookrightarrow\Theta$ with $\Phi$ and $\Theta$ finite
and $b\in\Phi$. Let $B$ be an open ball (not necessarily an element of $\Theta$)
such that $B\cap b\not=\emptyset.$

Then, we can add a ball $B^*\subset B$ to $\Theta$ and a branch $b^*\dashv^* b$ to $\Phi$ in such a way that:\\\\
{\bf (add1)} $B^*\not\in \Theta$;\\\\
{\bf (add2)} $b^*\not\in\Phi$;\\\\
{\bf (add3)} $\Phi\cup\{b^*\}\hookrightarrow\Theta\cup\{B^*\}$, with $\iota\(b^\ast\)=B^\ast$.
\end{lemma}
\proof
Consider
$$A=B\setminus\bigcup_{\cbra B'\in\Theta:\overline{B'}\subset\iota(b)\cket}\overline{B'}.$$
This set is open.
Now suppose $\overline{B'}\cap b\not=\emptyset$ for some $B'\subset \iota(b)$. As we observed in the caption to Figure \ref{tangent}, this implies that $\partial B'$ is tangent to $b$ (because $B'\cap b=\emptyset$). Therefore, $\overline{B'}\cap b=\partial{B'}\cap b$, which consists of a single point. Since there are only finitely many elements in $\Theta,$
not all points of $b\cap B$ (which is infinite) lie on some $\overline{B'}$, hence $b\cap A\not=\emptyset$ and we can fix ${\bf y}\in b\cap A$. 

Let $\delta>0$ be small enough so that
$$B_\delta({\bf y})\subset A.$$

Then it is evident that $b^*$ and $B^*$ can be picked inside $B_\delta({\bf y})$, satisfying all the conditions we wanted; just take $B^*$ so that $\partial B^*$ is tangent to $b$, and let $b^*\subset B^*$ be the radius that touches $b$.
\endproof

Suppose now that we have
$$\xymatrix{ \Theta^0\ar[r]^-{f} & \Phi^{1}\ar@{^{(}->}[r]&\Theta^{1}},$$
so that $f_-$ maps each element of $\Theta^0$ to its corresponding branch in $\Phi^1$ and there is a map $\iota$ embedding $\Phi^1$ into $\Theta^1$. Composing these two we get a map
$$f_\circ:\Theta^0\to\Theta^1,$$
which preserves inclusion; $f_\circ(B)$ is the smallest element of $\Theta^1$ containing $f(B)$.

The mapping $\iota$ does not need to be surjective; in fact, since $\bigcup_{b\in \Phi}b$ is
nowhere dense, we can add open balls to $ \Theta^1$ using $\varepsilon$-refinements; we will exploit this fact later on to create an open map $\omega$ into our Kripke frame ${\bf K}=\left <W,R,g\right>$. Further, as the lemma below shows, the resulting map $\dot\omega$ will satisfy
$$\dot\omega f=g\dot \omega,$$
so we can apply Theorem \ref{prel}. This is the strategy we will follow throughout the rest of this paper.
\begin{lemma}\label{commutes}
Let
$${\bf K}=\left <W,R,g\right>$$
be a dynamic Kripke frame with limits such that
$$g\left (\lim_{n\to\infty}w_n\right )=\lim_{n\to\infty}g\left (w_n\right ).$$
Let
$$\xymatrix{ \Theta^0\ar[r]^-{f} & \Phi^{1}\ar@{^{(}->}[r]&\Theta^{1}}$$
and $\omega:\Theta^0\cup\Theta^1\to W$ be continuous on each $\Theta^n$
and satisfy $g\omega=\omega f_\circ$.

Then, if ${\bf x}\in root(\Theta^0),$
$$g\dot\omega({\bf x})=\dot\omega f({\bf x}).$$
\end{lemma}
\proof
Let ${\bf x}\in root(\Theta^0)$. Recall that
$$root(\Theta^0)=dom(\Theta^0)\cup lim(\Theta^0);$$
we will split the proof into these two cases.\\\\
Case 1: ${\bf x}\in dom(\Theta^0)$.\\

Putting together Definitions \ref{treemaps} and \ref{embedded}, we have the following commutative diagram:
\begin{equation}
\xymatrix{
\Theta^0\ar@{-->}[r]^-{f_-} & \Phi^1\ar@{-->}[r]^-{\iota} & \Theta^1\\
dom(\Theta^0 )\ar[u]^-{\theta^0}\ar[r]^-{f} & im(\Phi^1)\ar[u]^{\phi^1}\ar@{^{(}->}[r]
& dom(\Theta^1)\ar[u]^-{\theta^1}
.}
\end{equation}

Then, by diagram chasing, along with the assumption that $g\omega=\omega f_\circ$, we can see that
$$
\begin{array}{lllll}
\dot\omega(f({\bf x}))&=&\omega(\theta^1(f({\bf x})))&&\\\\
&=&\omega(f_\circ(\theta^0({\bf x})))&&\\\\
&=&g\left (\omega(\theta^0({\bf x}))\right )&&\\\\
&=&g(\dot\omega({\bf x})).&&\\
\end{array}
$$
This is what we wanted.\\\\
Case 2: ${\bf x}$ is a limit point of $\Theta^0$.\\

We claim that
$$\theta^1_n(f({\bf x}))=f_\circ(\theta^0_n({\bf x})).$$

To prove this, note that $dom\(\Theta\)$ is dense in $root\(\Theta\)$, and
$$f\(\theta^0_n({\bf x})\cap dom\(\Theta\)\)\subset f_\circ(\theta^0_n({\bf x})).$$
Because $f$ is continuous,
$$f({\bf x})\in \overline{f_\circ(\theta^0_n({\bf x}))}$$
for all $n$. But
$$\overline{f_\circ\(\theta^0_{n+1}({\bf x})\)}\subset f_\circ\(\theta^0_n({\bf x})\),$$
hence $f({\bf x})\in f_\circ(\theta^0_n({\bf x}))$ for all $n$.

Now we must check that all elements of $\Theta^1$ containing ${\bf x}$ are of this form. Assume
$$B=\theta^1_m\(f({\bf x})\)$$
for some $m$.

Take $N$ large enough so that
$$f_\circ\(\theta^0_N({\bf x})\)\subset B;$$
such an $N$ exists because $diam\(\theta^0_n({\bf x})\)$ tends to zero as $n$ tends to infinity.

Then, by the last line of Definition \ref{embedded}, there exists some branch $b\in\Phi^1$ such that $\iota(b)=B$ and
$$f_-(\theta^0_N({\bf x}))\dashv^* b.$$
Since $f_-$ is onto, we then get
$$b=f_-(\theta^0_n({\bf x}))$$
for some $n\leq N$, hence
$$B=f_\circ(\theta^0_n(\bf x)).$$
Thus for all $m$, $\theta^1_m\(f({\bf x})\)$ is of the desired form.

We conclude that
$$f\(\theta^0_n\({\bf x}\)\)=\theta^1_n\(f\({\bf x}\)\)$$
for all $n$.

Using this and the fact that $g$ commutes with limits, we see that\\\\
$$
\begin{array}{lllll}
\dot\omega(f({\bf x}))&=&\lim_{n\to\infty}\omega(\theta^1_n(f({\bf x})))&&\\\\
&=&\lim_{n\to\infty}\omega(f_\circ(\theta^0_n({\bf x})))&&\\\\
&=&\lim_{n\to\infty}g\left (\omega(\theta^0_n({\bf x}))\right )&&\\\\
&=&g\left (\lim_{n\to\infty}\omega(\theta^0_n({\bf x}))\right )&&\\\\
&=&g(\dot\omega({\bf x})).&&\\
\end{array}
$$
\endproof
\section{Dynamic sequences of open ball trees and completeness results for $\R^2$}

Throughout this section, ${\bf K}=\left<W,R,g\right>$ denotes a dynamic Kripke frame where limits are defined,
$g$ commutes with limits and every world has finitely many $R$-successors.

We are now interested in constructing a dynamic analogue of open ball trees.
\begin{definition}[dynamic sequence]
Let
$$\Omega =\left \{\{\Theta^n\}_{n=0}^\infty,\{\Phi^n\}_{n=1}^\infty,f,\omega\right \},$$
where $\{\Theta^n\}_{n=0}^\infty$ is a sequence of open ball trees,
$\{\Phi^n\}_{n=1}^\infty$ a sequence of segment trees, $$f:\overline{\bigcup_{n=0}^\infty root(\Theta^n)}\to \R^2$$
a continuous function and $$\omega:\bigcup_{n=0}^\infty\Theta^n\to W$$ continuous on each
$\Theta^n$ and satisfying $g\omega=\omega f$.

Assume $root(\Theta^n)$ has radius 1 for all $n$ and $n\not=n'$ implies
$$\overline{root(\Theta^n)}\cap\overline{root(\Theta^{n'})}=\emptyset.$$

If
$$\xymatrix{ \Theta^n\ar[r]^-{f} & \Phi^{n+1}\ar@{^{(}->}[r]&\Theta^{n+1}},$$
we say $\Omega$ is a {\em dynamic sequence} of open ball trees.
\end{definition}
 
Note that Lemma \ref{commutes} guarantees that $g\dot\omega=\dot\omega f$.
\begin{figure}[htp]
\begin{center}
\scalebox
{0.9}
{
\includegraphics[viewport=205 270 600 355]{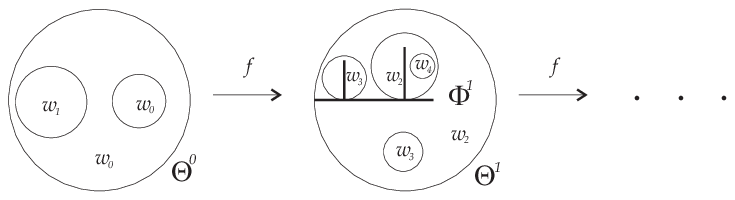}
}
\end{center}
\caption
{
The beginning portion of a dynamic sequence of open ball trees corresponding to the Kripke frame on Figure \ref{dynk}.
}
\end{figure}

We will use the following abuse of notation: if $B=B_{\varepsilon}({\bf x})$, we will write
$2B=B_{2\varepsilon}({\bf x})$.\\\\
\begin{definition}[$\varepsilon$-extendibility]\label{ten}
Let
$$\Omega =\left \{\{\Theta^n\}_{n=0}^\infty,\{\Phi^n\}_{n=1}^\infty,f,\omega\right \}$$
be a dynamic sequence of open ball trees, where each $\Theta^n$ is finite. Let $\tilde\Theta^n$ be an open ball tree
and $\tilde\omega:\tilde\Theta^n\to W$ a continuous map.

We say that $\(\tilde \Theta^n,\tilde \omega\)$ is {\em $\varepsilon$-extendible} over
$(\Theta^n,\omega)$
if\\\\
{\bf (ext1)} if $B,B'\in\tilde\Theta^n\setminus\Theta^n$, $2B\cap2B'\not=\emptyset$ implies
$B=B'$;\\\\
{\bf (ext2)} if $n>0$, $\Phi^n\hookrightarrow\tilde\Theta^n$;\\\\
{\bf (ext3)} if $B=B_\delta({\bf x})\in\tilde\Theta^n\setminus\Theta^n$, then
$$f(B_{2\delta}({\bf x}))\subset B_\eta(f({\bf x})),$$
where
$$\eta=\frac{\varepsilon}{3}d\left (f({\bf x}),\partial \theta^{n+1}(f({\bf x}))\right ).$$

We say that $\(\tilde \Theta^n,\tilde \omega\)$ is an {\em $\varepsilon$-extendible refinement} of
$\(\Theta^n,\omega\)$
if it is both $\varepsilon$-extendible over and an $\varepsilon$-refinement of $(\Theta^n,\omega)$.\\
\end{definition}

Roughly, $\varepsilon$-extendibility allows us to add elements to the open ball trees of our sequence while making only small alterations of $f$. Condition ${\bf (ext1)}$ will allow us to realize several small deformations of $f$ independently of each other, while ${\bf (ext3)}$ will allow us to take limits of tree maps.

\begin{lemma}\label{prev}
Let
$$\Omega =\left \{\{\Theta^n\}_{n=0}^\infty,\{\Phi^n\}_{n=1}^\infty,f,\omega\right \}$$
be a dynamic sequence of open ball trees, $N\in\N$ fixed, and $\varepsilon>0$.

Then, there exists an $\varepsilon$-extendible refinement $\(\tilde\Theta^{N},\tilde\omega\)$ of
$\(\Theta^{N},\omega\)$.
\end{lemma}
\proof
The construction is very similar to the one we used to generate $\varepsilon$-refinements in the static case, but
requires a bit more care.

Let $E=im(\Phi^n)$ if $n>0$, $E=\emptyset$ otherwise. Note that $E$ is closed and contains no open balls.

Let $C_1,...,C_I$ be the open, connected sets of the $\varepsilon$-grid $\Gamma_\varepsilon\(\(\Theta^N\)^\circ\setminus E\).$
As before, $\dot\omega$ is equal to a constant $\omega(C_i)$ on each $C_i$.

Let $w_{i,0},...,w_{i,M_i}$ be all the successors
of $\omega(C_i)$, and pick out $M_i$ disjoint open balls $B_{i,m}$ 
contained in $C_i$.

Now, suppose $B_{i,m}=B_{\rho_{i,m}}({\bf x}_{i,m})$, and let
$$\eta_{i,m}=\frac{\varepsilon}{3}d\left (f({\bf x}_{i,m}),\partial \theta^{N+1}(f({\bf x}_{i,m}))\right ).$$

Because $f$ is continuous, we can find
$\delta_{i,m}\in\(0,\frac{\rho_{i,m}}{2}\)$ such that
$f(B_{2\delta_{i,m}}({\bf x}_{i,m}))\subset B_{\eta_{i,m}}(f({\bf x}_{i,m}))$.

Then set
$$\tilde\Theta^N=\Theta^N\cup\bigcup_{i,m}\{B_{\delta_{i,m}}({\bf x}_{i,m})\}$$
and
$$\check\omega(B_{\delta_{i,m}}({\bf x}))=w_{i,m}.$$

It is evident that $\(\tilde\Theta^N,\tilde\omega\)$ satisfies {\bf (ext1)} and {\bf (ext3)}. As for {\bf (ext2)}, suppose $N>0$.
Note that $B\in\tilde\Theta^N\setminus\Theta^N$ implies $B\cap im(\Theta^n)=\emptyset$ by construction; hence
if
$\xymatrix{ \Phi^{n+1}\ar@{^{(}->}[r]^-{\iota}&\Theta^{n+1}},$ it still holds that
$\xymatrix{ \Phi^{n+1}\ar@{^{(}->}[r]^-{\iota}&\tilde\Theta^{n+1}}.$ 
\endproof
\begin{definition}[$\varepsilon$-extensions of tree maps]
With the notation of Definition \ref{ten}, let
$$\tilde f:\overline{root(\Theta^n)}\to\R^2$$
be continuous and $\Xi^n\supset\Theta^n$, $\Xi^{n+1}\supset\Theta^{n+1}$
and $\Psi^{n+1}\supset\Phi^{n+1}$ be such that $\Psi^{n+1}\hookrightarrow\Xi^{n+1}$.

We say
$$\xymatrix{ \Xi^n\ar[r]^-{\tilde f} & \Psi^{n+1}\ar@{^{(}->}[r]&\Xi^{n+1}}$$
is an $\varepsilon$-extension of
$$\xymatrix{ \Theta^n\ar[r]^-{f} & \Phi^{n+1}\ar@{^{(}->}[r]&\Theta^{n+1}}$$
if\\\\
{\bf (etm1)} ${\tilde f_-}|_{\Theta^n}=f_-$;\\\\ 
{\bf (etm2)} $\left \|\tilde f-f\right \|_\infty<\varepsilon$;\\\\
{\bf (etm3)} 
$$d\(f({\bf x}),\partial \theta^{n+1}(f({\bf x}))\)>(1-\varepsilon)d\(\tilde f({\bf x}),\partial \theta^{n+1}(f({\bf x}))\).$$\\
\end{definition}
\begin{lemma}\label{next}
Let
$$\Omega =\left \{\{\Theta^n\}_{j=0}^\infty,\{\Phi^n\}_{n=1}^\infty,f,\omega\right \}$$
be a dynamic sequence of open ball trees, $N\in\N$ and $\varepsilon>0$.

Assume $\left (\tilde\Theta^N,\tilde\omega^N\right )$ is $\varepsilon$-extendible over
$\left (\Theta^N,\omega\right )$.

Then there exist a function $\tilde f$, a segment tree $\tilde\Phi^{N+1}\supset\Phi^{N+1}$
and $\(\tilde\Theta^{N+1},\tilde\omega^{N+1}\)$ $\varepsilon$-extendible over 
$\(\Theta^{N+1},\omega\)$ such that
$$\xymatrix{ \Theta^N\ar[r]^-{f} & \Phi^{N+1}\ar@{^{(}->}[r]&\Theta^{N+1}}$$
is an $\varepsilon$-extension of
$$\xymatrix{ \tilde\Theta^N\ar[r]^-{\tilde f} & \tilde\Phi^{N+1}\ar@{^{(}->}[r]&\tilde\Theta^{N+1}}.$$
\end{lemma}
\proof
Let
$$\tilde\Theta^N\setminus\Theta^N=\{B_{\delta_1}({\bf x}_1),...,B_{\delta_M}({\bf x}_M)\}.$$

Because $\left (\tilde\Theta^N,\tilde\omega^N\right )$ is $\varepsilon$-extendible over
$\left (\Theta^N,\omega\right )$,
$$f\(B_{2\delta_m}({\bf x}_m)\)\subset B_{r_m}(f({\bf x}_m)),$$
where
$$r_m=\frac{\varepsilon}{3}d\left (f({\bf x}_m),\partial \theta^{N+1}(f({\bf x}_m))\right ).$$

Let $\rho_m\in(0,r_m)$ be small enough so that
$$f\left (B_{\rho_m}(f({\bf x}_m)))\subset B_{\eta_m}(f^2({\bf x})\right ),$$
where
$$\eta_m=\frac{\varepsilon}{3}d\left (f^2({\bf x}_m),\partial \theta^{N+2}(f^2({\bf x}_m))\right ).$$
By successive iterations of Lemma \ref{addbranch}, find balls
$B_m^*\subset B_{\rho_m}(f({\bf x}_m))$
and segments $b_m^*\subset B_m^*$ such that
$$\Phi^{N+1}\cup\bigcup_{m=1}^M\{b_m^*\}\hookrightarrow\Theta^{N+1}\cup\bigcup_{m=1}^M\{B_m^*\}.$$

If necessary, replace these by segments $\tilde b_m\subset b_m^*$ and open balls $\tilde B_m\subset B_m^*$
small enough so that $2\tilde B_m\cap 2\tilde B_{m'}\not=\emptyset$ implies $m=m'$.
This can be done because there are finitely many balls $B_m^*$.

Set
$$\tilde\Theta^{N+1}=\Theta^{N+1}\cup\bigcup_{m=1}^M\cbra\tilde B_m\cket$$
and
$$\tilde\Phi^{N+1}=\Phi^{N+1}\cup\bigcup_{m=1}^M\cbra\tilde b_m\cket.$$

Defining $\tilde\omega^{N+1}\left (\tilde B_m\right )=
g\left (\omega\left (B_{\delta_m}\left ({\bf x}_m\right )\right )\right),$ it is clear that
$\(\tilde\Theta^{N+1},\tilde\omega\)$ is $\varepsilon$-extendible over $\(\Theta^{N+1},\omega\)$ by the same arguments as used in 
Lemma \ref{prev}.

Because the balls $B_{2\delta_m}\({\bf x}_m\)$ are all disjoint, we can use Lemma \ref{pinch} $M$ times to obtain $\tilde f$ such that
$$\tilde f\succ_{B_{2\delta_m}({\bf x}_m)}^{B_{\delta_m}({\bf x}_m)\to\tilde b_m}f$$
on $\overline{B_{2\delta_m}({\bf x}_m)}$,
$$\left\|\tilde f|_{B_{2\delta_m}({\bf x}_m)}-f|_{B_{2\delta_m}({\bf x}_m)}\right\|_{\infty}<\varepsilon
d\left (f({\bf x}_m),\partial \theta^{N+1}(f({\bf x}_m))\right ),$$
and $\tilde f=f$ everywhere else.

The function $\tilde f$ thus defined satisfies {\bf (etm1)} by the definition of pinching,
and
$$\xymatrix{ \tilde\Theta^n\ar[r]^-{\tilde f} & \tilde\Phi^{n+1}\ar@{^{(}->}[r]&\tilde\Theta^{n+1}}.$$

Now, since $$d\left (f({\bf x}_m),\partial \theta^{N+1}(f({\bf x}_m))\right )\leq 1,$$
we have that
$$\left\|\tilde f-f\right\|_{\infty}<\varepsilon,$$
so {\bf (etm2)} holds.

Also note that, if ${\bf y}\in B_{\delta_m}({\bf x}_m)$, $\theta^N(f({\bf y}))=\theta^N(f({\bf x}_m))$
so
$$
\begin{array}{rl}
d\(\tilde f({\bf y}),\partial \theta^N(f({\bf y}))\)
>&d\(f({\bf x}_m),\partial \theta^N(f({\bf x}_m))\)-r_m\\\\
=&d\(f({\bf x}_m),\partial \theta^N(f({\bf x}_m))\)-\frac{\varepsilon}{3}
d\(f({\bf x}_m),\partial \theta^N(f({\bf x}_m))\)\\\\
=&\left (1-\frac{\varepsilon}{3}\right )d\(f({\bf x}_m),\partial \theta^N(f({\bf x}_m))\)\\\\
>&(1-\frac23\varepsilon-\frac{\varepsilon^2}{3})d(f({\bf x}_m),\partial \theta^N(f({\bf x}_m)))\\\\
=&(1-\varepsilon)\left (1+\frac\varepsilon 3\right )d(f({\bf x}_m),\partial \theta^N(f({\bf x}_m)))\\\\
=&(1-\varepsilon)\left (d(f({\bf x}_m),\partial \theta^N(f({\bf x}_m)))+r_m\right )\\\\
>&(1-\varepsilon)d(f({\bf y}),\partial \theta^N(f({\bf y}))).\\\\
\end{array}$$

Hence, {\bf (etm3)} holds as well and
$$\xymatrix{ \Theta^N\ar[r]^-{f} & \Phi^{N+1}\ar@{^{(}->}[r]&\Theta^{N+1}}$$
is an $\varepsilon$-extension of
$$\xymatrix{ \tilde\Theta^N\ar[r]^-{\tilde f} & \tilde\Phi^{N+1}\ar@{^{(}->}[r]&\tilde\Theta^{N+1}}.$$

\endproof

\begin{lemma}\label{partial}
Let $\Omega$ be as before, and $\varepsilon>0$.

Then, for all $N\in\N$, there exist sequences $\cbra\check\Theta^n\cket_{n=0}^{N+1}$ and
$\cbra\check\Phi^n\cket_{n=1}^{N+1}$, along with continuous maps
$$\check f:\bigcup_{n=0}^N\overline{root(\Theta^N)}\to\R^2$$
and
$$\check\omega:\bigcup_{n=0}^{N+1} \Theta^n\to W$$
such that\\\\
{\bf (lem1)} $\check\omega|_{\Theta^n}=\omega|_{\Theta^n}$ for $n\leq N$, and
$\check\Theta^n$ is $\frac{\varepsilon}{2}$-dense in $\tilde\Theta^n$;\\\\
{\bf (lem2)} $$\xymatrix{ \check\Theta^n\ar[r]^-{\check f} & \check\Phi^{n+1}\ar@{^{(}->}[r]&\check\Theta^{n+1}}$$
is an $\varepsilon$-extension of 
$$\xymatrix{ \Theta^n\ar[r]^-{f} & \Phi^{n+1}\ar@{^{(}->}[r]&\Theta^{n+1}}$$
for $n<N$;\\\\
{\bf (lem3)} $\(\check\Theta^{N+1},\check\omega|_{\check\Theta^{N+1}}\)$ is $\frac{\varepsilon}{2}$-extendible over
$\(\Theta^{N+1},\check\omega|_{\Theta^{N+1}}\)$.
\end{lemma}
\proof
We work by induction on $N$.

For the base case, use Lemma \ref{prev}
to find a $\frac{\varepsilon}{2}$-extendible refinement $\(\check\Theta^0,\check\omega|_{\check\Theta^0}\)$ of
$(\Theta^0,\omega).$

Then use Lemma \ref{next} to define a function
$$\check f:\overline{root(\Theta^0)}\to\R^2$$
such that
$$\xymatrix{ \check\Theta^0\ar[r]^-{\check f} & \check\Phi^{1}\ar@{^{(}->}[r]&\check\Theta^{1}}$$
is an $\varepsilon$-extension of
$$\xymatrix{ \Theta^0\ar[r]^-{f} & \Phi^{1}\ar@{^{(}->}[r]&\Theta^{1}},$$
where $\(\check\Theta^1,\check\omega|_{\check\Theta^1}\)$ is $\frac{\varepsilon}{2}$-extendible over
$(\Theta^1,\omega|_{\Theta^1})$.

For the induction step, assume that $\cbra\tilde\Theta^n\cket_{n=0}^{N+1}$, $\cbra\tilde\Phi^n\cket_{n=1}^{N+1}$,
$\tilde f$ and $\tilde\omega$ satisfy {\bf (lem1)}-{\bf (lem3)}.

Use Lemma \ref{next} to find extensions
$$\xymatrix{ \tilde\Theta^n\ar[r]^-{\tilde f} & \tilde\Phi^{n+1}\ar@{^{(}->}[r]&\tilde\Theta^{n+1}}$$
of
$$\xymatrix{ \Theta^n\ar[r]^-{f} & \Phi^{n+1}\ar@{^{(}->}[r]&\Theta^{n+1}}$$
such that $\tilde\Theta^{N+2}$ is $\frac{\varepsilon}{2}$-extendible over $\Theta^{N+2}$.

Use Lemma \ref{prev} to find an $\frac{\varepsilon}{2}$-extendible refinement
$\left (\check\Theta^{N+1},\check\omega|_{\check\Theta^{N+1}}\right )$
of $\left (\tilde\Theta^{N+1},\tilde\omega|_{\tilde\Theta^{N+1}}\right )$.

Once again use Lemma \ref{next} to find an extension
$$\xymatrix{ \check\Theta^n\ar[r]^-{\check f} & \check\Phi^{n+1}\ar@{^{(}->}[r]&\check\Theta^{n+1}}$$
of
$$\xymatrix{ \tilde\Theta^n\ar[r]^-{\tilde f} & \tilde\Phi^{n+1}\ar@{^{(}->}[r]&\tilde\Theta^{n+1}}$$
where $\check\Theta^{N+2}$ is $\frac{\varepsilon}{2}$-extendible over $\tilde\Theta^{N+2}$.

It is not hard to see that
$$\xymatrix{ \check\Theta^n\ar[r]^-{\check f} & \check\Phi^{n+1}\ar@{^{(}->}[r]&\check\Theta^{n+1}}$$
is an $\varepsilon$-extension of
$$\xymatrix{ \Theta^n\ar[r]^-{f} & \Phi^{n+1}\ar@{^{(}->}[r]&\Theta^{n+1}}.$$

Setting $\check f|_{\overline{root(\Theta^{n})}}=\tilde f|_{\overline{root(\Theta^{n})}}$
and $\check\Theta^n=\tilde\Theta^n$ for $n<N$, we see that
$\cbra\check\Theta^n\cket_{n=0}^{N+2}$, $\cbra\check\Phi^n\cket_{n=1}^{N+2}$,
$\check f$ and $\check\omega$ satisfy {\bf (lem1)}-{\bf (lem3)}.
\endproof

\begin{corollary}\label{cor}

Let
$$\Omega=\left \{\{\Theta^n\}_{n=0}^\infty,\{\Phi^n\}_{n=1}^\infty,f,\omega\right \}$$
be a dynamic sequence of finite open ball trees.

Then, there exists another dynamic sequence of finite open ball trees
$$\check\Omega=\left \{\{\check\Theta^n\}_{n=0}^\infty,\{\check\Phi^n\}_{n=1}^\infty,\check f,\check\omega\right \}$$
such that for all $n\geq 0$ the following conditions hold:\\\\
{\bf (cor1)} $\Theta^n\subset\check\Theta^n,$ $\check\omega|_{\Theta^n}=\omega|_{\Theta^n}$ and
$\check\Theta^{n}$ is $\varepsilon$-dense in $\Theta^n$;\\\\
{\bf (cor2)}
$$\xymatrix{ \check\Theta^n\ar[r]^-{\check f} & \check\Phi^{n+1}\ar@{^{(}->}[r]&\check\Theta^{n+1}}$$
is an $\varepsilon$-extension of
$$\xymatrix{ \Theta^n\ar[r]^-{f} & \Phi^{n+1}\ar@{^{(}->}[r]&\Theta^{n+1}}.$$
\end{corollary}
\proof
Since we can construct, by Lemma \ref{partial}, partial sequences satisfying {\bf (lem1)}-{\bf (lem3)} for all $N$, it
is easy to see that we can take the union of all these sequences to obtain $\check\Omega$.
\endproof

If $\check\Omega$ and $\Omega$ are as in the corollary, we will say $\check\Omega$ is {\em $\varepsilon$-dense} in $\Omega$.

\begin{lemma}\label{isdyn}
Let
$$\Omega^m =\left \{\{\Theta^{n,m}\}_{n=0}^\infty,\{\Phi^{n,m}\}_{n=1}^\infty,f^m,\omega|_
{\bigcup_{n=0}^\infty\Theta^{n,m}}\right \}$$
be a sequence of dynamic sequences of finite open ball trees such that $\Omega^{m+1}$
is $\left (\frac{1}{2}\right )^{m+1}$-dense in $\Omega^m$.

Let
$$\Theta^n=\bigcup_{m=0}^{\infty}\Theta^{n,m},
\Phi^n=\bigcup_{m=1}^\infty\Phi^{n,m} \text{ and }
f=\lim_{m\to\infty}f^m.$$
Then,
$$\Omega =\left \{\{\Theta^{n}\}_{n=0}^\infty,\{\Phi^n\}_{n=1}^\infty,f,\omega\right \}$$
is a dynamic sequence of open ball trees and $\dot\omega$ is continuous and open.
\end{lemma}
\proof
Openness and continuity of $\dot\omega$ is given by Lemma \ref{old}.

The main thing we must establish is that $f$ is a tree map.
Clearly it is continuous, since the sequence $\left \{f^m\right \}_{m=0}^\infty$
is Cauchy by {\bf (etm2)}.

Let ${\bf x}\in dom(\Theta^n)$.
Pick $M$ large enough
so that $B=\theta^n({\bf x})\in\Theta^{n,M}$.

We first show that
$$\theta^{n+1}(f({\bf x}))=f^{M}_\circ(B)$$
(recall that $f^{M}_\circ(B)$ is the smallest element of $\Theta^{n+1,M}$ containing $f^M(B)$).

It is clear that
$$f({\bf x})\in \overline{f^{M}_\circ(B)},$$
so we only need to prove that
$$d_m=d\(f^m\({\bf x}\),\partial f^{M}_\circ\(B\)\)$$
does not tend to zero as $m$ tends to infinity.

If it did, there would be some $m_0$ such that $d_{m_0}=\sup_m d_m$.

But then, by {\bf (etm3)},
$$d_{m+1}>\left (1-\left ( \frac{1}{2}\right )^{m+1}\right )d_m$$
for all $m$, so that 
$$d_m-d_{m+1}<\left ( \frac{1}{2}\right )^{m+1}d_m\leq\left ( \frac{1}{2}\right )^{m+1}d_{m_0}$$
and hence
$$
\begin{array}{lcl}
\lim_{m\to\infty}d_m&=&d_{m_0}-\sum_{m=m_0}^\infty \(d_{m}-d_{m+1}\)\\\\
&>&d_{m_0}-\sum_{m=m_0}^\infty\left ( \frac{1}{2}\right )^{m+1}d_{m_0}\\\\
&\geq&0;\\\\
\end{array}
$$
this contradicts the assumption.

Therefore
$$\lim_{m\to\infty}f^m({\bf x})\in f_\circ(B),$$
which is what we wanted.

On the other hand, for all $m>M$,
$$f^m({\bf x})\in {f^M_-}(B),$$ and because ${f^M_-}(B)$
is closed,
$$f({\bf x})\in f_-(\theta^n({\bf x}))={f^M_-}(B).$$

We conclude that
$$\phi^{n+1}(f({\bf x}))= f_-(\theta^n({\bf x})),$$
and hence $f$ is a tree map.

The rest of the properties that $\Omega$ must satisfy are trivially verified.
\endproof

We finally have all the tools we need to prove our main result.
\begin{theorem}
Let $\varphi$ be any ${\bf DTL}$ formula. If $\varphi$ can be satisfied in a locally finite Kripke frame, then $\varphi$ can be satisfied on $\R^2$.
\end{theorem}
\proof
Let $\varphi$ be any ${\bf DTL}$ formula and suppose $M=\left <W,R,g,V\right >$ is a locally finite dynamic Kripke frame where 
$$w_0\models_M\varphi.$$

We can suppose without loss that $W$ has limits that commute with $g$.\\

Let
$$\Theta^{n,0}=
\left \{B_1\left (\left [
\begin{array}{c}
3n\\
0\\
\end{array}
\right ]\right )\right \}.$$

For $n>0$, let $\Phi^{n,0}=\{ root(\Phi^{n,0}) \}$, where
$$root(\Phi^{n,0}) = \left \{
\left [
\begin{array}{c}
x\\
y\\
\end{array}
\right ]\in\R^2:
3n-1\leq x\leq 3n
,\text{ }
y=0 
\right \}.$$

Let
$$f^0
\bra
\begin{array}{c}
x\\
y\\
\end{array}
\ket
=
\left [
\begin{array}{c}
3(n+1)- \left \|  \left [
\begin{array}{c}
x\\
y\\
\end{array}
\right ]-
\left [
\begin{array}{c}
3n\\
0\\
\end{array}
\right ]\right\| \\
\\
0\\
\end{array}
\right ]
\text{ if }
\left[
\begin{array}{c}
x\\
y\\
\end{array}
\right ]\in root(\Theta^{n,0})$$
\\
(where $\|\cdot\|$ denotes the Euclidean norm) and\\
$$\omega(root(\Theta_0^n))=g^n(w_0).$$
\\\\

Let
$$\Omega^0 =\left \{\{\Theta^{n,0}\}_{n=0}^\infty,\{\Phi^{n,0}\}_{n=1}^\infty,f^0,\omega|_
{\bigcup_{n=0}^\infty\Theta^{n,0}}\right \}.$$

By successive iterations of Corollary \ref{cor}, we can construct a sequence $\{\Omega^m\}_{m=0}^\infty$
such that $\Omega^{m+1}$ is $\left (\frac 1 2\right )^{m+1}$-dense in $\Omega^m$.

By Lemma \ref{isdyn},
$$\Omega =\left \{\left \{\bigcup_{m=0}^{\infty}\Theta^{n,m}\right \}_{n=0}^\infty,
\left \{\bigcup_{m=1}^\infty\Phi^{n,m}\right \}_{n=1}^\infty,
\lim_{m\to\infty}f^m,\omega\right \}$$
is a dynamic sequence of open ball trees, and defines an open and continuous mapping
$$\dot\omega:\bigcup_{n=0}^\infty {root(\Theta^{n,0})}\to W$$
satisfying $g\dot\omega=\dot\omega f.$

The function $f$ can be extended to all of $\R^2$ using Tietze's extension theorem, giving us a dynamic model
$$M^{\R^2}=\left <\R^2,f,\omega^{-1}V\right >.$$

Pick ${\bf x}\in root(\Theta^{0,0})$ such that $\dot\omega({\bf x})=w_0$.\\

Then, by Theorem \ref{prel},
$${\bf x}\models_{M^{\R^2}}\varphi.$$
\endproof
\begin{corollary}
${\bf DTL^\bigcirc}$, ${\bf DTL_1}$, and all other fragments of ${\bf DTL}$ that are complete for locally finite Kripke frames are complete for $\R^2$.
\end{corollary}
\proof
Immediate from the previous theorem and Theorem \ref{aleks}.
\endproof
\section{Incompleteness of {\bf DTL} for $\R^2$}

In this section we present a formula of {\bf DTL} that is valid in all locally connected topological spaces (and hence in $\R^2$), but not in arbitrary spaces. This shows that the results we have proven so far cannot be extended to arbitrary formulas of {\bf DTL}.

Let
$$\varphi=\Box(\ast(\Box p\vee\Box\neg p))\wedge\ast p\rightarrow\Box \ast p.$$

We claim that if ${\bf M}=\<\<X,{\bf T}\>,f,V\>$ is a dynamical topological model and $\<X,{\bf T}\>$ is locally connected, then
$${\bf M}\models \varphi.$$

To see this, assume there is a point $x\in X$ such that
$$x\models \Box(\ast(\Box p\vee\Box\neg p))\wedge\ast p,$$
and let $U_x$ be a connected neighborhood of $x$ such that
$y\models \ast(\Box p\vee\Box\neg p)$
for all $y\in U_x$. Then for all $n$,
$$U_x\subset f^{-n}(V(p)^\circ)\cup f^{-n}(V(\neg p)^\circ).$$

These two sets are disjoint and open, so by connectedness of $U_x$ we conclude that
$$U_x\subset f^{-n}(V(p)^\circ)\text{ or }U_x\subset f^{-n}(V(\neg p)^\circ).$$
But $x\models \ast p$, so the second case is impossible and
$$U_x\subset f^{-n}(V(p)^\circ)$$
for all $n$, which in turn implies that
$$x\models \Box\ast p,$$
so
$$x\models\varphi.$$

However, a countermodel to $\varphi$ is not difficult to construct on a locally disconnected space. Let
$$Z=\cbra 0\cket\cup\cbra\frac 1n:n\in \N\cket$$
with the usual topology and $X=Z\times \N$.
Set
$$f(x,n)=(x,n+1)$$
for all $(x,n)\in X$ and
$$V(p)=\cbra(x,n):n=0\text{ or }x<\frac 1n\cket.$$

It is then easy to see that
$$(x,0)\not\models\varphi.$$

\appendix
\section{Notions from analysis}

The results and notions from analysis used in this paper are all fairly standard and can be found in a text such as \cite{rudin}. Here we give a quick summary of the necessary background and the conventions that we use.

We denote the usual Euclidean norm on $\R^2$ by $\left\|{\bf x}\right\|$, and $d$ stands for the metric induced by this norm, i.e.,
$$d({\bf x},{\bf y})=\left\|{\bf x}-{\bf y}\right\|.$$
If $Y\subset \R^2$, we set
$$d({\bf x},Y)=\inf\cbra d({\bf x},{\bf y}):{\bf y}\in Y\cket.$$

The diameter of a set $X$ is denoted by $diam(X)$; that is,
$$diam(X)=\sup\cbra d({\bf x},{\bf y}):{\bf x},{\bf y}\in X\cket.$$
In the case of a ball in $\R^2$ this is just the usual diameter.

If ${\bf x}\in\R^2$ and $\varepsilon>0$, then $B_\varepsilon({\bf x})$ denotes the open ball centered at ${\bf x}$ with radius $\varepsilon$.

If $C\subset \R^2$ is compact and $f:\R^2\to\R^2$ is a continuous function, the {\em infinity norm} of $f$ is defined by
$$\left\|f\right\|_\infty=\sup\cbra\left\|f({\bf x})\right\|:{\bf x}\in C\cket.$$
This defines a norm on the vector space of all such functions, and the set of continuous functions is a complete metric space under this norm; that is, if $\cbra f_n\cket_{n<\infty}$ is a Cauchy sequence, there exists a continuous function $f:C\to\R^2$ such that
$$\lim_{n\to\infty}f_n=f,$$
where the latter convergence holds pointwise. Recall that the sequence being Cauchy means that, for all $\varepsilon>0$, there exists $N>0$ such that $n,m>N$ implies that
$$\left\|f_n-f_m\right\|_\infty<\varepsilon.$$

We also need the following result, which we state here in the form that we use it:
\begin{theorem}[Tietze extension theorem]
If $\<X,d\>$ is a metric space, $A\subset X$ is closed and
$$f:A\to \R^n$$
is continuous, there exists a continuous map
$\hat f : X\to\R^n$
such that
$$\hat f|_{A}=f|_A.$$
\end{theorem}

\bibliographystyle{plain}

\end{document}